\documentclass[12pt]{article}
\usepackage{amssymb}
\usepackage{amsfonts}

\input{tcilatex}
\begin{document}

\begin{center}
\bigskip

\textbf{Some remarks on structural matrix rings }

\textbf{and matrices with ideal entries}
\end{center}

\bigskip

\begin{center}
Stephan Foldes

Tampere University of Technology, PL 553, 33101 Tampere, Finland

sf@tut.fi

\bigskip

Gerasimos Meletiou

TEI\ of Epirus, PO Box 110, 47100 Arta,\ Greece

gmelet@teiep.gr

\bigskip

December 2010\bigskip

\textbf{Abstract}
\end{center}

\textit{Associating to each pre-order on the indices }$1,...,n$\textit{\ the
corresponding structural matrix ring, or incidence algebra, embeds the
lattice of }$n$-\textit{element} \textit{pre-orders into the lattice of }$%
n\times n$ \textit{matrix rings. Rings within the order-convex hull of the
embedding, i.e. matrix rings that contain the ring of diagonal matrices, can
be viewed as incidence algebras of ideal-valued, generalized pre-order
relations. Certain conjugates of the upper or lower triangular matrix rings
correspond to the various linear orderings of the indices, and the incidence
algebras of partial orderings arise as intersections of such conjugate
matrix rings. }

\bigskip

Keywords: structural matrix ring, incidence algebra, pre-order, quasi-order,
triangular matrix, conjugation, semiring, ideal lattice, subring lattice%
\[
\]

\textbf{1 Conjugate subrings}

\bigskip

\textit{Rings }are understood to be possibly non-commutative, and to have a 
\textit{unit} (multiplicatively neutral) element, which is assumed to be
distinct from the zero (null) element and is denoted by $1$ or $I$ or a
similar symbol$.$ \textit{Subrings} are understood to contain the unit
element. An $n\times n$ \textit{matrix} is viewed as a map defined on the
set $\mathbf{n}^{2}=\left\{ 1,...,n\right\} ^{2}.$ (Here $n\geq 1$ is
assumed.) The ring of $n\times n$ matrices over a ring $R$ is denoted by $%
M_{n}(R).$

A \textit{pre-order}, also called \textit{quasi-order,} is a reflexive and
transitive binary relation $\lesssim $\ on a set $S$, an \textit{order} (or 
\textit{partial order}) is an anti-symmetric pre-order, and a \textit{linear
(}or\textit{\ total) order} is an order in which any two elements are
comparable. Instead of the generic notation $\lesssim ,$ specific pre-orders
may be denoted by other symbols such as $\theta $. Given a ring $R$ and a
pre-order $\lesssim $ on $\mathbf{n}$, the \textit{structural matrix ring }$%
M_{n}(\lesssim ,R)$ over $R$ is defined by 
\[
M_{n}(\lesssim ,R)=\left\{ \text{ }A\in M_{n}(R):\text{ }\forall i,j\text{ \
\ }A(i,j)=0\text{ unless }i\lesssim j\text{ }\right\} 
\]%
The full matrix ring $M_{n}(R)$, the subrings of all upper triangular
(respectively lower triangular), and of all diagonal matrices, are examples
of structural matrix rings. Structural matrix rings are essentially the same
as incidence algebras of finite pre-ordered sets, although the latter term
is sometimes used under the assumption that the base ring $R$ is a field [F]
or that the pre-order in question is an order, possibly on an infinite set,
as in [R]. Ring-theoretical properties of incidence algebras most relevant
to the present context were studied in [DW, F, MSW, W1, W2].

The set of all pre-orders defined on any given set constitutes a lattice
whose minimum is the equality relation.

\bigskip

\textbf{Proposition 1} \ \textit{For any ring} $R,$\textit{\ the map
associating to each pre-order }$\lesssim $\textit{\ on} $\mathbf{n}$ \textit{%
the corresponding structural matrix ring }$M_{n}(\lesssim ,R)$\textit{\
provides an embedding of the lattice of pre-orders on} $\mathbf{n}$ \textit{%
into the lattice of subrings of the matrix ring} $M_{n}(R).$ \textit{The
embedding preserves also infinite greatest lower and least upper bounds.\ \
\ \ \ }

\bigskip

\textbf{Proof \ }Preservation of greatest lower bounds is obvious. The
preservation of upper bounds is a consequence of the description of the
least upper bound of a family of pre-orders as the transitive closure of the
least binary relation that is implied by the pre-orders in the family. $\ \
\ \ \ \square $

\bigskip

The lattice embedding described in the proposition above is generally not
surjective, and for $n\geq 2$\ its range is order-convex if and only if $R$
is a division ring. Section 2 will provide a description of the order-convex
hull of the embedding's range.

Subrings $S$ and$\ T$ of $M_{n}(R)$\ are said to be \textit{permutation
conjugates }if there is a permutation matrix $P$ such that%
\[
S=\left\{ PAP^{-1}:A\in T\right\} 
\]%
Such subrings are obviously isomorphic under the automorphism $A\mapsto
PAP^{-1}$ of $M_{n}(R).$

All the $n!$ linear orders on the finite set $\mathbf{n}$ are isomorphic,
and the isomorphisms among them are precisely the self-bijections of the
underlying set $\mathbf{n.}$ Consequently we have:

\bigskip

\textbf{Proposition 2} \ \textit{For any pre-order} $\lesssim $ \textit{on} $%
\mathbf{n}$ \textit{and any ring} $R,$\ \textit{the following conditions are
equivalent}:

(i) $\lesssim $ \textit{is a linear order},

(ii) $M_{n}(\lesssim ,R)$ \textit{is a permutation conjugate of the ring of
upper triangular matrices,}

(iii) $M_{n}(\lesssim ,R)$ \textit{is a permutation conjugate of the ring of
lower triangular matrices. \ \ \ \ \ \ \ \ \ \ \ }$\square $

\bigskip

The well-known fact that every partial order is the intersection of its
linear extensions yields:

\bigskip

\textbf{Proposition 3 \ }\textit{For any pre-order} $\lesssim $ \textit{on} $%
\mathbf{n}$ \textit{and any ring} $R,$ \textit{the following conditions are
equivalent}:

(i) $\lesssim $ \textit{is an order},

(ii) $M_{n}(\lesssim ,R)$ \textit{is the intersection of some permutation
conjugates of the ring of upper triangular matrices,}

(iii) $M_{n}(\lesssim R)$ \textit{is the intersection of some permutation
conjugates of the ring of lower triangular matrices. \ \ \ \ }$\square $

\bigskip

These considerations were first developed, in the context of fields, in
[FM1]. They were related to a ring property concerning one-sided and
two-sided inverses in [FSW], addressing questions originating in [C] and
[SW].

\bigskip

\[
\]

\textbf{2 \ Matrices with ideal entries}

\bigskip

The order-convex hull of the embedding provided by Proposition 1 consists
obviously of those matrix rings that contain the ring of diagonal matrices.
We now show that these rings can be viewed as generalized incidence
algebras, corresponding to generalized relations that are the analogues of
pre-order relations in a reticulated semiring-valued framework.

\bigskip

By a \textit{semiring} we understand a set endowed with a binary operation
called \textit{sum} (denoted additively) and a binary operation called 
\textit{product} (denoted multiplicatively) such that:

(i) the sum operation defines a commutative semigroup,

(ii) the product operation defines a semigroup,

(iii) both distributivity laws $a(b+c)=ab+ac$ and $(b+c)a=ba+ca$ hold for
all members $a,b,c$ of the underlying set.

\bigskip

The set $\mathcal{I}(R)$\ of (bilateral) ideals of any ring\ $R$\ is a
semiring under ideal sum and product of ideals. This semiring is
lattice-ordered by inclusion. The set $M_{n}(\mathcal{I}(R))$\ of $n\times n$
matrices with entries in $\mathcal{I}(R)$ is again a semiring under the
obvious sum and product operations, also lattice-ordered by entry-wise
inclusion:\ lattice join coincides with semiring sum, denoted $+$, while the
lattice meet operation $\wedge $ is entry-wise intersection$\mathbf{.}$ The
lattice $M_{n}(\mathcal{I}(R))$\ is complete, it is isomorphic to the $n^{2}$%
-th Cartesian power of the complete lattice $\mathcal{I}(R).$ The matrix $%
\mathbf{I,}$\ with the improper ideal $(1)=R$ in diagonal positions and the
trivial ideal $(0)$ in off-diagonal positions, is multiplicatively neutral
in the semiring $M_{n}(\mathcal{I}(R))$.

\bigskip

For every matrix $\mathbf{U}$ with ideal entries, $\mathbf{U}\in M_{n}(%
\mathcal{I}(R))$, consider the set of matrices%
\[
G=\left\{ A\in M_{n}(R):\text{ }\forall i,j\text{ \ \ }A(i,j)\in \mathbf{U}%
(i,j)\right\} 
\]%
This set is always an additive subgroup of $M_{n}(R),$ and it is a subring
if and only if $\ \mathbf{U}^{2}+\mathbf{I}\leq \mathbf{U}$ \ in the ordered
semiring of $n\times n$ matrices with ideal entries. In that case, we say
that the subring $G$ is \textit{defined by} $\mathbf{U.}$

\bigskip

There is a natural lattice embedding from the lattice of all pre-orders on $%
\mathbf{n}$ into the lattice $M_{n}(\mathcal{I}(R)),$ namely the map $\theta
\mapsto \mathbf{U}$ where $\mathbf{U(}i,j)$ is the improper or trivial ideal
of $R$ according to whether the relation $i\theta j$ holds or not. The
matrix $\mathbf{U}$\ so obtained always satisfies $\mathbf{U}^{2}+\mathbf{I}%
\leq \mathbf{U}$, and the subring of $M_{n}(R)$ defined by it coincides with
the structural matrix ring $M_{n}(\theta ,R).$ For this reason we denote by $%
M_{n}(R,\mathbf{U})$ the matrix ring defined by any $\mathbf{U}\in M_{n}(%
\mathcal{I}(R))$ satisfying $\mathbf{U}^{2}+\mathbf{I}\leq \mathbf{U}.$ Any $%
\mathbf{U}\in M_{n}(\mathcal{I}(R))$ is viewed as an $\mathcal{I}(R)$-valued
relation on the $n-$element set $\mathbf{n}$, the inequality $\mathbf{U}^{2}+%
\mathbf{I}\leq \mathbf{U}$ generalizes transitivity and reflexivity of
relations, and $M_{n}(R,\mathbf{U})$ may then be thought of as the incidence
algebra of the generalized pre-order $\mathbf{U}$.

Matrices with ideal entries $\mathbf{U}\in M_{n}(\mathcal{I}(R))$\
satisfying the generalized transitivity condition $\mathbf{U}^{2}\leq 
\mathbf{U}$\ only were used by van Wyk [W1], and recently again by Meyer,
Szigeti and van Wyk [MSW], in the description of ideals of stuctural matrix
rings .

Any matrix with ideal entries $\mathbf{U}\in M_{n}(\mathcal{I}(R))$
satisfying $\mathbf{U}^{2}+\mathbf{I}\leq \mathbf{U}$ is said to be \textit{%
reflexive-transitive. }If\textit{\ }$R$ is a division ring, then the range
of the map $\theta \mapsto \mathbf{U}$ described above, embedding the
pre-order lattice into $M_{n}(\mathcal{I}(R)),$ is exactly the set of
reflexive-transitive matrices.

For any ring $R,$\ the set of reflexive-transitive matrices with ideal
entries constitutes a complete lattice under the ordering of $M_{n}(I(R))$
(but not a sublattice of $M_{n}(I(R))$\ in general)$.$

\bigskip

\textbf{Proposition 4 \ }\textit{For any ring }$R$, \textit{the map }$%
\mathbf{U\mapsto }M_{n}(R,\mathbf{U})$ \textit{establishes a lattice
isomorphism between:}

(i) \textit{the lattice }$\left\{ \mathbf{U}\in M_{n}(\mathcal{I}(R)):%
\mathbf{U}^{2}+\mathbf{I}\leq \mathbf{U}\right\} $\textit{\ of }$n\times n$%
\textit{\ reflexive-transitive matrices with ideal entries, }

(ii) \textit{the lattice of subrings of }$M_{n}(R)$ \textit{containing all
diagonal matrices.}

\bigskip

\textbf{Proof} \ Obviously if $\mathbf{U}^{2}+\mathbf{I}\leq \mathbf{U}$\
then $M_{n}(R,\mathbf{U})$ contains all diagonal matrices, and if $\mathbf{%
U\subseteq V}$ then $M_{n}(R,\mathbf{U})\subseteq M_{n}(R,\mathbf{V}).$

Conversely, let $N\subseteq M_{n}(R)$ be any matrix ring including all
diagonal matrices. For every $1\leq i,j\leq n$ the set $U_{ij}=\left\{
A(i,j):A\in N\right\} $\ is an ideal of $R$, the matrix $\mathbf{U}=(U_{ij})$
with ideal entries can be seen to satisfy $\mathbf{U}^{2}+\mathbf{I}\leq 
\mathbf{U,}$ and $M_{n}(R,\mathbf{U})=N.$ \ \ $\ \ \ \ \ \ \ \ \square $

\bigskip

The proof above is presented in [FM2] in somewhat greater detail, together
with some consequences. For division rings, where there is only the trivial
and the improper ideal, the structure of the lattice\ of subrings of $%
M_{n}(R)$ containing all diagonal matrices is independent of the choice of
the particular division ring $R.$ In the general case, it is the ideal
lattice structure of $R$ that determines the structure of this upper section
of the subring lattice of $M_{n}(R).$

\bigskip 
\[
\]

\textbf{References}

\bigskip

[C] P.M. Cohn, Reversible rings, Bull. London Math. Soc. 31 (1999) 641--648

\bigskip

[DW] S. Dascalescu, L. van Wyk, Do isomorphic structural matrix rings have
isomorphic graphs? Proc. Amer. Math. Soc. 124 (1996) 1385--1391.

\bigskip

[F] R.B. Feinberg, Polynomial identities of incidence algebras, Proc. Amer.
Math. Soc. 55 (1976) 25--28.

\bigskip

[FM1] \ S. Foldes, G. Meletiou, On incidence algebras and triangular
matrices, Rutcor Res.Report 35-2002, Rutgers University, 2002. Available at
rutcor.rutgers.edu/\symbol{126}rrr

\bigskip

[FM2] \ S. Foldes, G. Meletiou, On matrix rings containing all diagonal
matrices, Tampere University of Technology, August 2007
http://math.tut.fi/algebra/

\bigskip

[FSW] \ S. Foldes, J. Szigeti, L. van Wyk, Invertibility and Dedekind
finiteness in structural matrix rings. Forthcoming in Linear and Multilin.
Algebra 2011. Manuscript iFirst 2010, 1--7,
http://dx.doi.org/10.1080/03081080903357653

\bigskip

[MSW] J. Meyer, J. Szigeti, L. van Wyk: On ideals of triangular matrix
rings, Periodica Math. Hung. Vol. 59 (1) (2009), 109-115

\bigskip

[R] \ G.-C. Rota, On the foundations of combinatorial theory I. Theory of M%
\"{o}bius functions, Zeitschrift Wahrscheinlichkeitstheorie 2 (1964) 340--368

\bigskip

[SW] J. Szigeti, L. van Wyk, Subrings which are closed with respect to
taking the inverse, J. Algebra 318 (2007) 1068--1076

\bigskip

[W1] L. van Wyk, Special radicals in structural matrix rings, Communications
Alg. 16 (1988) 421-435

\bigskip

[W2] L. van Wyk, Matrix rings satisfying column sum conditions versus
structural matrix rings, Linear Algebra Appl. 249 (1996) 15--28

\end{document}